\documentclass[10pt]{amsart}

\usepackage{amsmath}
\usepackage{graphicx}
\usepackage{amsthm}
\usepackage{graphicx}
\usepackage{verbatim}
\usepackage{enumerate}
\usepackage{mathtools}

\theoremstyle{plain}
\newtheorem{thm}{Theorem}[section]

\newtheorem{lem}[thm]{Lemma}
\newtheorem{pr}[thm]{Proposition}
\newtheorem{cor}[thm]{Corollary}
\newtheorem{defn}[thm]{Definition}

\newtheorem{conj}[thm]{Conjecture}

\newtheorem*{theorem*}{Theorem}
\theoremstyle{remark}

\newtheorem*{remark*}{Remark}

\font\elevenss=cmss11

\font\eightss=cmss8

\font\sixss=cmss8 at 6pt

\newfam\ssfam
\textfont\ssfam=\elevenss \scriptfont\ssfam=\eightss
 \scriptscriptfont\ssfam=\sixss

\def\Z{\mathbb{Z}}

\def\R{\mathbb{R}}

\def\ee{\varepsilon}
\def\E{{\mathbb E}}
\def\P{{\mathbb P}}
\def\Cox{\hfill \Box}
\def\disp{\displaystyle}

\def\|{{\, | \, }}

\def\pois{{\mathcal P}}
\def\poi{{\mathcal M}}

\def\simp{\bigtriangleup}
\def\vv{{\bf v}}
\def\xx{{\bf x}}
\def\yy{{\bf y}}
\def\nn{{\bf n}}
\def\sss{{\bf s}}
\def\ttt{{\bf t}}
\def\cL{{\mathcal L}}

\def\CM{{\mathcal M}}
\def\TM{\Omega}
\def\TMm{\TM_0}
\def\TMo{\TM^o}

\begin{document}

\title{Coarsening in One Dimension: \\ Invariant and asymptotic states}

\author{Emanuel A. Lazar}
\address{Laboratory for Research on the Structure of Matter, 
University of Pennsylvania, 
Philadelphia, PA 19104}
\email{mlazar@seas.upenn.edu}
\urladdr{www.seas.upenn.edu/~mlazar/}

\author{Robin Pemantle}
\address{Department of Mathematics, University of Pennsylvania, 
Philadelphia, PA 19104}
\email{pemantle@math.upenn.edu}
\urladdr{www.math.upenn.edu/~pemantle/}

\begin{abstract}
We study a coarsening process of one-dimensional cell complexes.  
We show that if cell boundaries move with velocities proportional 
to the difference in size of neighboring cells, then the average cell
size grows at a prescribed exponential rate and the Poisson distribution 
is precisely invariant for the distribution of the whole process, rescaled
in space by its average growth rate.  We present numerical
evidence toward the following universality conjecture: starting 
from any finite mean stationary renewal process, the system when 
rescaled by $e^{-2t}$ converges to a Poisson point process.
For a limited case, this makes precise what has been observed 
previously in experiments and simulations, and lays the foundation 
for a theory of universal asymptotic states of dynamical cell complexes.
\end{abstract}

\maketitle

\vspace{-0.2in}

\setcounter{equation}{0}
\section{Introduction}
\label{sec:intro}

Many physical systems can be abstracted as cell complexes whose geometry 
and topology change over time through deterministic evolution equations.  
Soap foams and polycrystalline metals are two examples of such systems 
that evolve under a generalized curvature flow to reduce an energy 
associated with the co-dimension 1 cells.  Over time, some cells grow 
while others shrink and disappear.  This is sometimes called a 
coarsening process, because as cells disappear, the average 
size of the remaining cells monotonically increases~\cite{lazar2011evolution}.  

A remarkable {\bf self-similarity} property has long been observed in 
both experiments and numerical simulations of these systems: if an 
arbitrary cell complex is allowed to evolve by fixed evolution 
equations, then its scale-independent properties will converge, 
in a statistical sense, to a time-independent set of properties.  
For example, the distribution of normalized cell sizes converges to 
some fixed distribution, even as the system itself continues to 
evolve \cite{lazar2011evolution, lazar2011more, mason2012statistical}.  
This general behavior is sometimes referred to as statistical 
self-similarity \cite{mullins1986statistical, mullins1989self}.  
A second remarkable {\bf universality} property is also observed: 
these self-similar properties are largely independent of the initial 
conditions, and instead depend on the deterministic evolution dynamics.  

\begin{figure}
\begin{center}    
\includegraphics[width=0.7\textwidth]{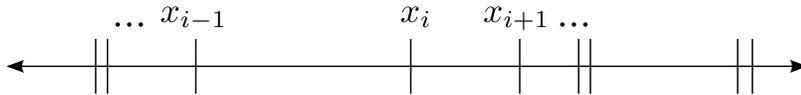}
\caption{A one-dimensional cell complex with a few of the cell boundaries 
$x_i$ labeled.}
\label{onedb}
\end{center}
\end{figure}

The primary aim of this paper is to give a rigorous analysis of a 
one-dimensional coarsening process in which a limit theorem can
be proved.  The second goal is to explore the degree to which the
limiting behavior is independent of initial conditions.  These will 
be carried out in the context of a particular one-dimensional system 
which can be viewed as one-dimensional curvature flow.  

We parametrize configurations by the locations $x_i$ of the boundaries 
between the 1-cells.  The particular system we initially consider is one 
in which cell boundaries move with a velocity proportional to the 
difference of the sizes of two neighboring 1-cells.  We let 
$X = (x_i) \in \mathbb{R}$ be an ordered set of points on a 
one-dimensional manifold, as illustrated in Figure \ref{onedb}.  
The evolution of this system is governed by two rules:
\begin{eqnarray}
\frac{dx_i}{dt} & = & -x_{i-1} + 2x_i - x_{x+1} \label{eq:rule1} \\[2ex]
\mbox{If } x_i = x_j \mbox{ for } i \neq j: && \hspace{-0.2in}
   \mbox{ we remove $x_j$ and reindex the remaining points.} \label{eq:rule2}
\end{eqnarray}
Rule 2 prevents points from crossing one another.  It is straightforward 
to see that 1-cells that are larger than their two neighbors grow, while 
those smaller than their neighbors shrink.   

\subsection*{Previous Work} 
Important work has been done in studying one-dimensional dynamical 
cell complexes, though these systems are often known by the names 
clustering or coagulation systems.  Carr and Pego~\cite{carr1992self} 
considered a model in which at discrete time steps the shortest interval 
in a finite system is taken and joined with its neighbors.  
They use Laplace transforms to establish the existence of a 
self-similar solution.  This system makes use of discrete time 
steps and ``long-range'' interaction. 
Derrida~\cite{derrida1995coarsening, derrida1997non} considered 
one-dimensional Ising and Potts models at zero temperature and considered 
the rates at which these systems coarsen.  He showed that if the initial 
condition is random, the system enters a scaling regime where the average 
size of domains grows as $t^{1/2}$.  Nontrivial exponents that govern 
other properties of these coarsening systems are also calculated.  
Hunderi and Ryum~\cite{hunderi1996influence, hunderi1979computer} 
considered a system in which $\disp \frac{dx_i}{dt} = \frac{1}{x_{i+1}-x_i} - 
\frac{1}{x_i-x_{i-1}}$.  They used numerical simulations to demonstrate 
that a size correlation develops between neighboring cells;  
this correlation calls into question the reliability of mean field 
theories \cite{fortes1992grain}, which tend to ignore size correlations 
between neighboring cells, in analyzing these systems.  

Mullins~\cite{mullins1991one} considered a slightly more 
general version of the system considered in this paper, and 
established some very preliminary results for these systems.  
Unfortunately, despite work of Mullins and of others, Mullins' 
conclusion, ``Very little is known, however, about the general 
conditions under which self-similarity is to be expected'' 
remains true nearly twenty-five years later.

\setcounter{equation}{0}
\section{Results}

Our main result is the invariance of the Poisson distribution
for cell boundaries under the evolution~\eqref{eq:rule1}--\eqref{eq:rule2}.  
We begin with an analysis of what happens with a finite number of 
such points, moving in a finite interval.  For boundary conditions,
here and throughout, we impose two extra particles that remain 
stationary at the endpoints~0 and~1.

\begin{thm} \label{th:finite}
Begin with $n$ points chosen uniformly on $[0,1]$.
The law after time $t$ is that of $Z$ points chosen 
uniformly and independently on $[0,1]$ where $Z$ is 
a binomial with parameters $n$ and $e^{-2t}$.
\end{thm}

Because these points are moving around, colliding with each other 
and coalescing, we refer to them henceforth as {\em particles}.
Some consequences follow immediately from Theorem~\ref{th:finite},
properties of the binomial and Poisson distributions, and 
scaling invariance. 
Let $\pois(\lambda , L)$ denote the law of a Poisson process of
intensity $\lambda$ times Lebesgue measure on $[-L,L]$; alternatively
it can be described as the law of $Z$ IID uniform points in $[-L,L]$
when $Z$ has Poisson distribution with mean $\lambda$.

\begin{cor} \label{cor:consequences}
~~\\[-5ex] \begin{enumerate}[1.]
\item Starting from $n$ particles uniformly distributed in $[-L,L]$,
the number of particles in a given subset $A \subseteq [-L,L]$ 
after time $t$ has binomial distribution
$${\rm Bin}\, \left ( n , e^{-2t} \frac{|A|}{2L} \right ) \, .$$
\item Starting from $\pois(\lambda , L)$, the law of the 
configuration after time $t$ is a Poisson configuration
with intensity $\lambda e^{-2t}$ times Lebesgue measure
on $[-L,L]$.
\item For $M \geq L > 0$, let $\cL (M , L)$ be the law of the 
time $t$ configuration restricted to $[-L , L]$, starting from 
$\pois(\lambda , M)$.  Then the law $\cL (M, L)$ does not depend on $M$.
\end{enumerate}
\end{cor}

Next we turn to the infinite process.  The first thing is to establish an
existence result.  This is not trivial because of the possibility 
of influence coming in from infinity in finite time and because
at collision times, which are dense, the trajectory of every particle
become non-analytic.
In Section~\ref{sec:poisson} we will define the space $\Omega$
of infinite particle configurations and prove the following result. 

\begin{thm}[existence of a weak solution for the infinite system] 
\label{th:existence}
Let $\TM$ be the space of continuous trajectories of infinite 
configurations defined in Section~\ref{sec:poisson}.  For every
$\lambda > 0$ there exists a weak solution to~\eqref{eq:rule1}
with initial condition Poisson of intensity $\lambda$, in the
sense of Definition~\ref{def:solution}.  Any such solution has the 
following properties.
\begin{enumerate}[(i)]
\item The trajectory of each particle is always differentiable 
from the right.
\item Excluding a measure zero set of initial configurations, the 
collision times will be distinct.  
\item At all times that a particle is not involved in a collision, 
that particle's position has a two-sided derivative.  
\end{enumerate}
\end{thm}

The infinite model allows for a statement of Poisson invariance 
cleaner than the one in Corollary~\ref{cor:consequences}.

\begin{thm} \label{th:poisson invariant}
Under any weak solution $\P_\lambda$, the time $t$ law of the 
configuration is Poisson with intensity $\lambda e^{-2t}$.
\end{thm}

\setcounter{equation}{0}
\section{Evolution on a finite interval}
\label{sec:finite}

In order to define the process on a finite interval, we need to
choose how the first and last particles evolve.  As previously
mentioned, rather than using periodic boundary conditions, we 
impose extra particles at the endpoints that never move.

The space of configurations of $n$ particles on $[0,1]$ is the simplex 
$$\Delta_n := \{ \xx = (0 , x_1 , \ldots , x_n , 1) : 
   0 = x_0 \leq x_1 \leq \cdots \leq x_n \leq x_{n+1} = 1 \} \, .$$
We include the degenerate simplex $\Delta_0 := \{ (0,1) \}$.
The interior $\Delta_n^o$ of $\Delta_n$ is the set of configurations with 
all particles distinct from each other and from the endpoints $\{ 0 , 1 \}$.
After a collision, two particles, say $x_i$ and $x_{i+1}$, get
stuck together.  The configuration is no longer in $\Delta_n^o$,
rather it is now on the boundary $\partial \Delta_n$.  The boundary
is composed of faces of various dimensions defined below as
images of open simplices under various embeddings.  On each of these,
the evolution satisfies a different rule.  Rather than think of this 
as a single discontinuous rule on the closed simplex, we think of the 
particle as entering the $(n-1)$-simplex at the first collision time.  

To make this precise, we define a set of embeddings as follows.
For integers $n \geq 0$ and $k \in [0,n+1]$, define the map
$\iota_{n,k} : \Delta_n \to \Delta_{n+1}$ by the formula
$$ \iota_{n,k}(0,x_1, \ldots , x_n,1) = \begin{dcases*}
        (0,0, x_1, \ldots, x_n,1)  & $k = 0$\\
        (0,x_1, \ldots, x_n, 1, 1)  & $k = n+1$\\
        (0,x_1, \ldots, x_k, x_k, \ldots, x_n,1) & $1 \leq k \leq n$
        \end{dcases*} $$
Define the projection map from $\Delta_n$ to the union of $\Delta_k^o$
for $k \leq n$ by letting $\pi (0, x_1 , \ldots x_n, 1)$ be the result of
omitting repeated entries.  
Then $\pi \circ \iota_{n-1,k} (\xx) = \xx$ for each $\xx \in \Delta_{n-1}^o$
and each $k \in [0,n]$.  More generally, any composition of maps
$\iota_{n-1,k_{n-1}} \circ \cdots \circ \iota_{m,k_m}$ is inverted 
by $\pi$ on $\Delta_m^o$.

We now define an evolution on the configuration space $\Delta := 
\biguplus_n \Delta_n^o$, the disjoint union of $n$-point 
configuration spaces on $[0,1]$.  The evolution is defined as
a set of time $t$ maps $\{ \Phi_t \}$ on $\Delta$ and constructed
by induction on the dimension $n$ of the stratum in $\Delta$.  
The base step of the induction is to define $\Phi_t (\xx) \equiv 0$
on $\Delta_0$.  

Let $A_n$ denote the $(n+2) \times (n+2)$ matrix defined by
$$A_n = \left [ \begin{array}{ccccccccccc} 
0 & 0 & 0 & 0 & \cdots & 0 \\
-1 & 2 & -1 & 0 & \cdots & 0 \\
0 & -1 & 2 & -1 & \cdots & 0 \\
\vdots & \vdots &&&& \vdots \\
0 & \cdots & 0 & -1 & 2 & -1 \\
0 & \cdots & 0 & 0 & 0 & 0 
\end{array} \right ] \, .$$
The differential equation $\xx' = A_n \xx$ defines a flow
on $\Delta_n^o$.  For $\xx \in \Delta_n^o$, let
$$\tau (\xx) := \inf \left \{ t : e^{t A_n} \xx \in \partial \Delta_n 
   \right \} \, .$$
We may now recursively define 
$$\Phi_t (\xx) := \left \{ \begin{array}{lr}
   e^{t A_n} \xx & t < \tau (\xx) \\[2ex]
   \Phi_{t - \tau (\xx)} \left ( \pi ( e^{\tau A_n} \xx ) \right )
   & t \geq \tau (\xx)
   \end{array} \right.  $$
Informally, flow by $\xx' = A_n \xx$ until you reach the boundary,
then collapse whichever points need to be collapsed and continue
inductively (run $\xx' = A_m \xx$ where $m$ is the new number of 
particles, and so on).  For all but a set of measure zero of 
initial points, each collision will reduce the number of particles
by precisely~1; however, the flow is well defined even for simultaneous
coalescences, in which case $\pi$ drops dimension by more than~1.

\begin{remark*}
One can define the evolution directly on $\Delta_n$. 
At time $\tau (\xx)$, let $F$ be the unique face of 
$\Delta_n$ to which $\xx$ is interior.  Let $\iota$
be the composition of maps $\iota_{m,k}$ inverted
by $\pi$ on $F$; the map $\iota$ is unique even though
it may be represented by a sequence of compositions in 
many ways.  Instead of mapping $\xx$ to the projection
$\pi (\Phi_{\tau (\xx)} (\xx))$ of its time $\tau$ image, 
continue the evolution via
$$\Phi_{\tau + s} (\xx) := \iota \left ( \Phi_s (\pi (e^{\tau A_n} \xx 
   ) \right ) \, .$$
Informally, the evolution on $\Delta_n$ keeps track not only of the
coalesced system but of which original particles coalesced at what
is now particle $i$.
\end{remark*}

\setcounter{equation}{0}
\section{The time $t$ pre-image of a finite configuration}

The forward flow goes until it hits the boundary, jumps
down a dimension, and continues.  The reverse flow goes
inward away from the boundary, but unlike the forward flow,
it is capable at any time of reinterpreting its position
in $\Delta_m$ as a position on $\partial \Delta_{m+1}$, 
in which case it jumps up a dimension and begins to flow 
inward from the boundary of the bigger simplex.  Not only
that, but it can choose to jump up in any of $m+2$ ways:
any of the particles can split into two, or a particle
can enter at~0 or~1.  In this section we characterize all
possible reverse trajectories.  First we check that the
flow, when not splitting, does indeed flow inward.

\begin{lem} \label{lem:within}
For $t > 0$, the map $\xx \mapsto e^{-t A_n} \xx$ is one to one 
on $\R^{n+2}$ and maps $\simp_n$ into its own interior.
\end{lem}

\noindent{\sc Proof:} Letting $y_j = x_{j+1} - x_j$, the ODE
$\xx' =  - A_n \xx$ induces an ODE $\yy' = B_n \yy$ where
$B_n$ is the $(n+1) \times (n+1)$ matrix with ones on the
first super- and sub-diagonal, and $-2$ on the main diagonal 
except that $B_{1,1} = B_{n+1,n+1} = -1$, not $-2$.
Each $(B_n \yy)_j$ is a linear combination
of coordinates of $\yy$ with the only negative contribution
coming from $y_j$; hence if $\yy \in (\R^+)^n$ and $y_j = 0$
then $(B_n \yy)_j \geq 0$.  By Nagumo's Theorem~\cite{Aubi1977}
(or see~\cite[Theorem~1]{hart1972}), solutions to 
$\yy' = B_n \yy$ never leave the nonnegative orthant. 
It follows that the coordinates $x_j$ never collide.
$\Cox$

Now fix integers $n \geq m \geq 0$ and a point $\yy \in \Delta_m^o$, 
and define the following notation for the time $t$ dimension $n$ 
pre-image:
$$\Phi_t^{-1} (\yy ; n) := 
   \{ \xx \in \Delta_n : \yy = \Phi_t (\xx) \} \, .$$
We show that $\Phi_t^{-1} (\yy ; n)$ has a one to one parametrization
by the sequence of times at which the split occurs and the choice
of which particle splits at each time.

To define this formally, let $T (t;k)$ be the $k$-simplex 
of vectors $\ttt$ satisfying $0 \leq t_1 \leq \cdots \leq t_k \leq t$.  
Define 
$$S(n;m) := \{ (k_1 , \ldots , k_{n-m}) : 0 \leq k_i \leq n-i+1 \} \, .$$
As time passes backward through $t_k$, one of the existing $n-k$ 
particles splits; these choices are what is encoded by an element of 
$S(n;m)$.

\begin{pr} \label{pr:pre-image}
Fix integers $m \leq n$ and a positive real $t$.  Let $k := n-m$.
There is a map $\eta$ from $\Delta_m \times T(t;k) \times S(n;k)$ 
to $\Delta_n$ such that if $\xx \in \eta (\yy , \ttt , \sss)$ then 
\begin{enumerate}[(i)]
\item $\Phi_t (\xx) = \yy$
\item For each $1 \leq j \leq k$, the first entry time $\tau_j$
of $\Phi_\cdot (\xx)$ into $\Delta_{n-j}$ is equal to $t_j$
\item For each $1 \leq j \leq k$ the coalescence at time $\tau_j$
occurs between coordinates $s_j$ and $s_{j-1}$.  If $\tau_j = \cdots 
= \tau_{j+r}$, this is taken to mean that the $r+1$ coalescences
at this time are those identifying $s_{j+i}$ with $s_{j+i-1}$
for $i=0, \ldots , r$.
\end{enumerate}
\end{pr}

\noindent{\sc Proof:} 
We construct the map $\eta$ explicitly.  Run the clock backwards 
starting at time $t$, evolving backwards via $e^{-t A_m}$.  
By Lemma~\ref{lem:within}, the point stays in $\Delta_m^o$.
When the clock reaches $t_k$, jump from
$\Delta_m$ to $\Delta_{m+1}$ via the embedding $\iota_{m,s_k}$.
Continue evolving backwards via $e^{-t A_{m+1}}$ until
the clock reaches time $t_{k-1}$, then jump up to 
$\Delta_{m+2}$ via $\iota_{m+1,s_{k-1}}$.  Continue
in this manner, reaching $\Delta_n$ at time $t_1$ and
evolving in $\Delta_n$ by $e^{-t A_n}$ for a backward
time $t_1$ until the clock says zero. 

With this construction of $\eta$,
the consequences of the proposition are easily checked by
induction on $k$.  When $k=0$ the vectors $\ttt$ and $\sss$
are empty and the only fact to check is that $e^{-tA}$
inverts $e^{tA}$.  For $k > 1$, evolving back in time
by $t - t_k$ yields a point $\yy' := e^{-t A_m} \yy$ and
a point $\yy'' := \iota_{m+1,s_k} \yy'$.  Apply the induction
hypothesis to $\xx = \eta (\yy', \ttt', \sss')$ where $\ttt'$
and $\sss'$ drop that last coordinates of $\ttt$ and $\sss$
respectively to see that $\Phi_{t_k-} (\xx) := \lim_{s \uparrow t_k}
\Phi_s (\xx) = \yy''$.  Thus,
$\Phi_t (\xx) = \Phi_{t-t_k} (\Phi_{t_k} (\xx)) = \Phi_{t-t_k} (\yy')
= \yy$.  One final coalescence occurs at time $t_k$ in coordinates
$s_k$ and $s_k - 1$.  This completes the induction.
$\Cox$

\setcounter{equation}{0}
\section{Invariance theorem for finite configurations}
\label{sec:uniform}

This section is devoted to the proof of Theorem~\ref{th:finite}
and its consequences.  We reduce it to Lemma~\ref{lem:density}
below, prove the further Lemma~\ref{lem:borel}, then prove 
Lemma~\ref{lem:density}.  
It will be helpful to define a certain continuous time Markov chain 
on the nonnegative integers.  It is a pure death chain, that is, 
transitions from $n$ are allowed only to $n-1$.  The rate at which 
$n$ transitions to $n-1$ is $2n$.  If $\{ X_t : t \geq 0 \}$ is such 
a Markov chain, one interpretation is that $X_t$ is the number 
of particles alive at time $t$, where each particle independently
dies at rate $2$ (meaning, after an exponential random time 
with mean $1/2$).  Denote the transition probabilities for
this chain by $p_t$.  Thus, $p_t (n,m)$
is the probability, starting with $n$ particles, that precisely
$m$ are alive at time $t$.  Saying it another way, 
\begin{equation} \label{eq:binom}
p_t (n,\cdot) \sim {\rm Bin}\,(n,e^{-2 t}) \, .
\end{equation}
Let $g(t_1, \ldots , t_k) = g (n;\, t_1 , \ldots , t_k)$ 
denote the density for the first $k$ transitions starting from $n$.  
An explicit formula for $g$ is
\begin{equation} \label{eq:g}
g (n;\, t_1 , \ldots , t_k) = 2^k n(n-1) \cdots (n-k+1) 
   \exp \left ( - \sum_{j=1}^k 2 (n-j+1) (t_j - t_{j-1}) \right )
\end{equation}
where $t_0 = 0$ by convention.

Normalized Lebesgue measure on $\Delta_n$ is the probability
measure $\mu_n$ whose density with respect to Lebesgue measure 
is $n!$.  We first prove a special case of Theorem~\ref{th:finite}.

\begin{lem} \label{lem:pois 1}
Let $B$ be any Borel subset of $\simp_n$ and let $t$ be any
positive real number.  Then
$$\mu_n \{ x \in \simp_n : \tau_\xx > t \mbox{ and }
   e^{tA} \xx \in B \} = e^{-2nt} \mu_n (B) \, .$$
Consequently, after time $t$ the probability of still
having $n$ points is $e^{-2nt}$, and conditional on 
this the points are independently and uniformly distributed.
\end{lem}

\noindent{\sc Proof:}  The divergence of the vector field $F(\xx) = A_n \xx$ 
on $\R^{n+2}$ is equal to $2n$.  It follows that the Jacobian of the map 
$e^{tA}$ on $\R^n$ is equal to $e^{2nt}$; this is easily deduced
from the fact that a divergence free flow has Jacobian equal to 
one~\cite[Proposition~18.18]{lee} applied to the map $e^{-2t} e^{tA}$.
Also, the time $t$ map on $\R^n$ is one to one and is inverted
by $e^{-tA}$~\cite[Theorem~17.8]{lee}.  It follows that
the Lebesgue measure of $e^{-tA} B$ is equal to $e^{-2nt}$ times
the Lebesgue measure of $B$.  By Lemma~\ref{lem:within}, the
set $e^{-tA} [B]$ is contained in $\simp_n$ and for each
$\xx \in e^{-tA} [B]$ the stopping time $\tau_\xx$ is greater 
than $t$ (apply the lemma for all $s \leq t$).  The conclusion
of the lemma follows.
$\Cox$

To see that that yields the special case of Theorem~\ref{th:finite}, 
let $B$ be a subset of the interior of $\Delta_n$.  Then the probability, 
starting with a uniform random point of $\Delta_n$, of being in $B$ at
time $t$ is the $\mu_n$ measure of points $\xx \in \Delta_n$ 
that evolve in time $t$ to a point in $B$.  By Lemma~\ref{lem:pois 1}
this is $e^{-2nt} \mu_n (B)$.  As the probability of
a ${\rm Bin}\,(n,e^{-2t})$ random variable being equal to $n$
is precisely $e^{-2nt}$, the conclusion of Theorem~\ref{th:finite}
is verified in this case.
$\Cox$

The full theorem is a consequence of the next result.  Extend
the notation for the map $\eta$ so that it acts on sets:
for $B \subseteq \Delta_m$, $T \subseteq T(t;k)$ and 
$S \subseteq S(n;k)$ we define 
$$\eta [B,T,S] := \{ \eta (\yy, \ttt, \sss) : \yy \in B, \ttt \in T ,
   \sss \in S \} \, .$$
In general, for $K \subseteq \Delta_m \times T(t;k) \times S(n;k)$ we 
denote $\eta [K] := \{ \eta (\yy, \ttt, \sss) : (\yy,\ttt,\sss) \in K \}$.

\begin{lem} \label{lem:density}
Let $m = n-k \leq n$, let $t > 0$, and let 
$K = K_0 \times S(n;k)$ be the product of a measurable subset of
$\Delta_m \times T(t;k)$ with all of $S(n;k)$.  Let $g$ 
be the density defined in equation~\eqref{eq:g}.  Then
\begin{equation} \label{eq:tbp2}
\mu_n (\eta(K)) = \int_K g(\ttt) \, d\mu_L (\yy) \, d\ttt  \, .
\end{equation}
\end{lem}

This lemma says that Lebesgue measure on $\Delta_n$
pulled back by $\eta$ yields the product of normalized
Lebesgue measure with the measure on $T(t;k)$ having density $g$,
provided that one sums over all possible embedding sequences $\sss$.
Before proving it, let us check that it implies 
Theorem~\ref{th:finite}.  Setting $K = B \times T(t;k)$ 
makes $\eta (K)$ the set of all $\xx$ for which $\Phi_t (\xx) \in B$.  
The integral in~\eqref{eq:tbp2} factors as a product in the 
$\yy$ and $\ttt$ variables, resulting in 
$$\mu_n (K) = \mu_m (B) p_t (n,m) \, .$$ 
The left-hand side is $\P (\Phi_t (\xx) \in B)$ and the 
right-hand side is the normalized Lebesgue measure of $B$
times the probability of ${\rm Bin}\,(n,e^{-2t}) = m$.
This proves Theorem~\ref{th:finite}.

There are two tricky points in establishing~\eqref{eq:tbp2}.  
One is seeing how the product of $m$ dimensional measure on 
$\Delta_m$ and $k$ dimensional measure on $T(t;k)$ maps to
$m+k = n$ dimensional measure on $\Delta_n$.  In fact the
bifurcation at a variable time $s \in [t , t+dt]$ results
in a product with a spatial interval proportional to the
normal velocity of of the flow toward the boundary of $\Delta_n$.
The second tricky point is why the total of these factors involving
normal velocities is always constant if one sums over all
sequences $\sss$.  For the first of these two tricky points 
we record the following lemma, for which all notation is 
declared to be local.

\begin{lem} \label{lem:borel}
Let $H$ be a subset of $\R^d$ which is locally a halfspace with
boundary $K$.  Let $B$ be a Borel subset of $K$, let $F$ be
a vector field on $H$ with $F \cdot \nn > 0$ on $K$ where $\nn$ 
is the inward pointing unit normal.  Let $\{ \Phi_t \}$ denote 
the flow associated with the differential equation $\xx' = F(\xx)$.  
Let $B_t \subseteq H$ be the set $\{ \Phi_s (\xx) : 
\xx \in B , 0 \leq s \leq t \}$.  Then
\begin{equation} \label{eq:normal}
\lim_{t \downarrow 0} \frac{|B_t|}{t} 
   = \int_B F(\xx) \cdot \nn (\xx) \, d\xx \, .
\end{equation}
\end{lem}

\noindent{\sc Proof:} The map $\Phi : K \times [0,t]$ taking 
$(\xx , s)$ to $\Phi_s (\xx)$ has Jacobian $J(\xx) = F(\xx) \cdot \nn (\xx)$
at the point $(\xx , 0)$ because it has a matrix representation 
in blocks of dimension $d-1$ and~1 of
$$\left [ \begin{tabular}{c|c} {\LARGE $I_{d-1}$} \hspace{0.2in} 
& * \\ \hline 0 & $F(\xx)$ 
   \end{tabular} \right ]$$
in local coordinates $K \times [0,\infty)$.  The measure of 
$B_t$ is equal to $\int_{B \times [0,t]} J(\xx,t) (d \xx \times dt)$.
Using continuity, $J(\xx,t) = (1+o(1)) J(\xx , 0)$ whence the
integral is $(t + o(t)) \int_B J(\xx,0) \, d\xx$ which matches
the right hand side of~\eqref{eq:normal}.
$\Cox$

Continuing the proof of Lemma~\ref{lem:density}, we 
observe that it suffices to show~\eqref{eq:tbp2} for 
rectangles $K = B \times T$.  Let $B$
be a small ball about a point $\yy$ and $T$ be a small rectangle 
$\{ c_j \leq t_j \leq c_j + \Delta : 1 \leq j \leq k \}$ about
a point $\ttt$.  This makes $F \cdot \nn$ roughly constant over 
the backward evolution of $B$ jumping at times in $T$, at any
fixed time $s \leq t$, as long as the embedding sequence 
$\sigma$ is held fixed.  Letting $t_{j+1} = t$
and $t_0 = 0$, we have
\begin{equation} \label{eq:normal2}
|\eta (B,T,\sigma)| \sim |B| |T| 
   \exp \left [ \sum_{j=0}^k - 2(n-j) (t_{j+1} - t_j) \right ]
   \prod_{j=1}^k \lambda_j \, .
\end{equation}
Here, $|\cdot|$ denotes Lebesgue measure, $F$ is the
vector field $F(\xx) = A_{n-j+1} \xx$, and $\lambda_j$ is 
$F \cdot \nn$ evaluated at the point where the transition is 
made from dimension $n-j+1$ to dimension $n-j$.  The asymptotic
equivalence is as the diameter of $B$ goes to zero.  The exponential 
reflects the fact that the backwards flow shrinks volume at 
a rate of $-2(n-j)$ over the time interval between $t_j$ and $t_{j+1}$.  

The last step is to sum over $\sigma$.  We do this inductively,
beginning with the sum over all values of $\sigma_1$.
Still supposing that $B$ is a small ball about $\yy$ and $T$
is a small rectangle about $\ttt$, we see that as time goes
backward from $t$ all the way to $t_1$, the backward evolution 
takes $\yy$ to a point $\yy' = (y_1' , \ldots y_{n-1}') \in 
\simp_{n-1}$.  The point $\yy'$ corresponds to one of $n+1$ 
possible points under the different possible choices of $\sigma_1 
= 0 , \ldots , n$, the $j^{th}$ of which is the point 
$(y_1' , \ldots , y_j' , y_j', \ldots , y_{n-1}')$ (where 
$j=0$ or $n$ corresponds to prepending~0 or appending~1).  
The identification $\sigma_j$ embeds $\simp_{n-1}$ in 
$\simp_n$ as the set of points whose $j$ and $j+1$ 
coordinates are equal.  The dual vector $\vv \mapsto \vv \cdot \nn$
is the functional $\vv \mapsto v_{j+1} - v_j$.  Therefore, at the point
$(y_1' , \ldots , y_j' , y_j', \ldots , y_{n-1}')$, the value
$\lambda_j = -F \cdot \nn$ is $(y_{j+1}' - 2 y_j' + y_j') - 
(y_j' - 2 y_j' + y_{j-1}')$ which simplifies to $y_{j+1}' - y_{j-1}'$.  
Note that we have used $-F$ because Lemma~\ref{lem:borel} refers to the
velocity of the backward flow.  To summarize, 
\begin{equation} \label{eq:gaps}
\mbox{The factor $\lambda_j$ is the sum of the two gaps on 
either side of the coalesced point.}  
\end{equation}
When $j=0$ or $n$ we obtain just one gap.  Summing over values 
of $\sigma_1$ from $0$ to $n$ gives twice the sum of all gaps, which
is just the constant~2.  

Inductively, still restricting $B$ to a small ball about $\yy$ and
$T$ to a small rectangle about $\ttt$, we see that the sum over 
$\sigma_1 , \ldots , \sigma_i$ of $\prod_{j=1}^i \lambda_j$ 
is $2^i$.  When $i$ reaches $k$, we see that the sum over all 
$\sigma$ of $\prod_{j=1}^k \lambda_j$ is $2^k$.  We have shown
that summing over $\sigma$ leads to a constant factor.  Letting
$\Sigma$ denote the set of all embedding sequences, we see that
that $\eta (\yy,\ttt,\Sigma) := \sum_{\sigma \in \Sigma} 
\eta (\yy, \ttt , \sigma)$ has density $2^k$
independent of $\yy$.  Integrating over $\yy$ is therefore trivial 
and leads, for $T$ a small rectangle, to
$$|\eta (B,T,\Sigma)| =  |T| \cdot |B| \cdot 2^k 
   \exp \left [ -\sum_{j=0}^k 2(n-j) (t_{j+1} - t_j) \right ] \, .$$
Recalling that $\mu_j$ is $j!$ times Lebesgue measure,
we may rewrite this as
\begin{eqnarray*}
\frac{1}{n!} \mu_n (\eta (B,T,\Sigma)) 
   & = & |T| \frac{\mu_m (B)}{m!} 2^k 
   \exp \left [ -\sum_{j=0}^k 2(n-j) (t_{j+1} - t_j) \right ]  \\[1ex]
\mbox{hence}&&\\
\mu_{n} (\eta (B,T,\Sigma)) 
& = & 2^k \frac{n!}{m!} |T| \mu_m (B) 
   \exp \left [ -\sum_{j=0}^k 2(n-j) (t_{j+1} - t_j) \right ]  \, .
\end{eqnarray*}
Comparing to~\eqref{eq:g} with $k=n-m$ shows that for small rectangles $T$,
$$\mu_n (\eta (B,T,\Sigma)) = g(\ttt) \; |T| \; \mu_m (B) \, .$$
Integrating over $T$ yields the result for general sets $T$,
which is~\eqref{eq:tbp2}, implying Lemma~\ref{lem:density} and 
Theorem~\ref{th:finite}.
$\Cox$

\setcounter{equation}{0}
\section{The time reversed Markov chain on finite configurations}
\label{sec:reverse}

The evolution of a finite configuration is deterministic.  Nevertheless,
for some choices of initial measure, projecting the configuration to
its cardinality produces a continuous time Markov chain.  The next
result is an immediate consequence of Corollary~\ref{cor:consequences}.
Recall that $\mu_n$ is the uniform measure on $\Delta_n$ and let 
$\poi_\lambda$ denote a mixture of laws $\mu_n$ when $n$ has
Poisson distribution with mean $\lambda$.

\begin{pr} \label{pr:markov}
For $\xx \in \Delta$, let $N(\xx,t)$ denote the number of particles
in the configuration at time $t$ started from $\xx$.  
\begin{enumerate}[(i)]
\item If $\xx$ is a random variable with law $\mu_n$, then 
$N(\xx , \cdot)$ is a pure death Markov chain on $\Z^+$ 
with initial law $\delta_n$ and jump rate at time $t$ 
from state $n$ given by $R(n,t) = 2n$ (hence transition
kernel given by~\eqref{eq:binom}).
\item If $\xx$ is a Poisson with law $\poi_\lambda$, 
then $N(\xx , \cdot)$ is a pure death Markov chain on $\Z^+$
whose initial distribution is Poisson$(\lambda)$ and whose
transition rate is again $2n$.
$\Cox$
\end{enumerate}
\end{pr}

We now describe the time reversals of these two chains.  Both
time reversals are pure birth chains.  The first chain is a 
time homogeneous chain.  It starts at~1 and has rate $q(k,t) := 2k$
of jumping from $k$ to $k_1$.  It is killed at the time $\sigma$ 
that it first reaches $n+1$.  
The second has jump rates that depend on time but not on state.
The initial distribution is Poisson with a specified mean 
$\beta > 0$.  The rate of an upward jump at time $s$ is
$r(k,s) := 2 \beta e^{2s} \, ds$.  

\begin{pr} \label{pr:reversals}
~
\begin{enumerate}[(i)]
\item Let $\nu_n$ be the law on right-continuous trajectories of 
the pure death Markov chain with death rate $2k$ from state $k$, 
started from $\delta_n$ and run until the time $\tau$ that it 
reaches zero.  For a trajectory $\omega$, denote the time reversal
by $s \mapsto \omega' (s) := \omega (\tau - s)$, made right continuous
by setting $\omega' (s) = \omega' (s^+)$ when $\tau - s$ is a 
jump time of $\omega$.  Then the $\nu_n$ law of $\omega'$ is the
same as the law of a pure birth Markov chain with birth rate
$2k$ from state $k$ and initial distribution $\delta_1$, on the time 
interval $[0,\sigma)$, where $\sigma$ is the hitting time on $n+1$.
\item Fix $t$ and $\lambda > 0$ and let $Q_\lambda$ be the law on
right-continuous trajectories of the pure death Markov chain with 
death rate $2k$ from state $k$, started from a Poisson of mean
$\lambda$ and run for time precisely $t$.  Let $\omega'$ again 
denote the right-continuous reversal of the trajectory $\omega$,
that is, $\omega' (s) = \omega (t-s)$ except at jump times.
Then the $Q_\lambda$ law of $\omega'$ is the same as the law of 
a pure birth Markov chain with initial distribution Poisson of 
mean $\beta := \lambda e^{-2t}$ and birth rate $2 \beta e^{2 \lambda s}$
at time $s$, run for time precisely $t$.
\end{enumerate}
\end{pr}

\noindent{\sc Proof:} Part~$(i)$: 
A trajectory $\omega$ is specified by its $n$ jump times, 
$s_1 , s_2 , \ldots, s_{n-1} , \tau$.  The density
of the trajectory with respect to $n$-dimensional Lebesgue measure is
$$2^n n! \exp \left ( - 2n s_1 - 2(n-1) (s_2 - s_1) - \cdots 
   - 2 (\tau - s_n) \right )  \, .$$
Similarly, parametrizing trajectories of the birth chain by the
jump times $(r_1, \ldots r_{n-1}, \sigma)$ gives a density of
$$2^n n! \exp \left ( -2 r_1 - 4 r_2 - \cdots - 2n (\sigma - r_{n-1})
   \right ) \, .$$
Setting $\sigma = \tau$ and $r_j = \tau - s_{n-j}$ for 
$1 \leq j \leq n-1$ reduces both exponents to 
$- 2 \tau - 2 \sum_{j=1}^{n-1} s_j$.

For part~$(ii)$, begin by observing that a death rate of $2n$
is equivalent to the $n$ particles each dying independently at
rate~2.  Thus $Q_\lambda$ is the law of Poisson-$\lambda$ many 
particles each with an independent death time whose law has
density $2e^{-2s}$ with an atom of size $e^{-2t}$ being still
alive at time $t$.  A trajectory may be specified by death times
together with the number remaining alive at time $t$, which is a 
Poisson process with intensity $2 \lambda e^{-2s} \, ds + \lambda 
e^{-2t} \delta_t$.  The birth chain is specified by these same 
parameters.  Clearly the number alive initially is Poisson with
the correct mean $\beta = \lambda e^{-2t}$.  The arrival process
has rate $2 \beta e^{2 \lambda s} = 2 \lambda e^{- 2 \lambda (t-s)}$
which agrees with the intensity of the point process of death times
at time $t-s$.  A Poisson process is completely specified 
by its intensity, establishing the distributional identity.
$\Cox$

Let $R(n,t)$ be a jump rate for a pure birth process.  We define 
an associated Markov transition kernel with state space $\Delta$
as follows.  First a sample path of the birth process is generated,
giving the jump times of the trajectory.  Between jumps, the
trajectory evolves deterministically via $\xx' = - A_n \xx$
when $\xx \in \Delta_n^o$.  At jump times, a uniform $[0,1]$
random variable $U$ is generated, and the particle closest to
$U$ splits (including the frozen particles at~0 and~1).  Let
$\P = \P_{\mu,\tau,R}$ denote the law of a trajectory with
jump times given by $R$, started from $\mu$, with each trajectory 
$\omega'$  stopped at time $\tau (\omega')$.

\begin{thm}[time reversal] \label{th:reversal}
~~\\[-4ex]
\begin{enumerate}[(i)]
\item Let $R_1 (k,t) := 2k$.  Define $\mu*$ to be the measure giving
probability $1/2$ to 
$\iota_{0,0} (0)$ and $1/2$ to $\iota_{0,1} (0)$.  In other words, 
one of the frozen particles in the empty configuration is chosen by 
fair coin-flip to duplicate.  Let $\tau$ be the hitting time
on $\Delta_{n+1}$.  Then the time reversal of $\P_{\mu*,\tau,R_1}$ 
is the law of forward evolutions started from $\mu_n$ and stopped 
upon hitting $\Delta_0$.
\item Let $R_2 (k,t) := 2\beta e^{2s}$.  Let $\tau$ be a fixed constant, 
$t$.  Then the time reversal of $\P_{\poi_\beta,\tau,R_2}$ is the law of 
forward evolutions started from $\poi_\lambda$ and stopped at time $t$.
\end{enumerate}
\end{thm}

\noindent{\sc Proof:}  Let $\Omega'$ denote the space of reversals
of finite trajectories of forward evolutions.  The time reversal of 
the forward evolution from $\mu_n$ stopped at $\tau$ (respectively 
the forward evolution from $\poi_\lambda$ stopped at $t$)
is a well defined measure on $M_1$ (respectively $M_2$) on 
$\Omega'$.  By Proposition~\ref{pr:reversals}, the splitting 
times for $M_1$ (respectively $M_2$) are given by the birth
process $R_1$ (respectively $R_2$).  The only thing left to
check is the law of the jump at the splitting times.

By~\eqref{eq:normal2}, the jumps are chosen independently and
proportionally to $\lambda_j$, the normal component of $F$ 
upon embedding by $\iota_j$.  By~\eqref{eq:gaps}, these are
proportional to the sums of consecutive pairs of gaps.  If
$0 = x_0 \leq x_1 \leq \cdots x_k \leq x_{k+1} = 1$ and 
$U$ is a uniform random variable on $[0,1]$ then the index 
of $j$ such that $x_j$ is the nearest point to $U$ has 
law proportional to the sums of pairs of consecutive gaps.
This completes the proof of the theorem.
$\Cox$

\setcounter{equation}{0}
\section{The infinite system and a weak solution}
\label{sec:poisson}

In this section we construct the infinite process, which is
a probability measure $\P$ on the space of trajectories of 
infinite point configurations on $\R$.  This requires constructing
the space $\CM$ of infinite point configurations, then the space
of trajectories $\TM$, then giving the probability measure $\P$.  
Although the usual choice for the space of infinite point 
configurations is that of counting measures, we use here
a more direct construction that lends itself to taking limits
of the finite configurations we have already defined.

Topologize $\Delta$ by treating the projection $\pi$ as describing
attachment maps from $\partial \Delta_n$ to $\bigcup_{k=0}^{n-1} 
\Delta_k$ for each $n$.
Notice that this identifies all $n$ faces of each $\partial \Delta_n$.
Also note that $\Delta$ is not compact: a sequence of
points in $\Delta_n^o$ for $n \to \infty$ has no limit.
In this topology the trajectories $\{ \Phi_t (\xx) \}$ are
continuous.  They are differentiable except at collision
times; one-sided derivatives exist even at collision times.

Let $\Delta_n^{[a,b]}$ denote the simplex of vectors of length $n$ 
of elements of $[a,b]$ whose coordinates are nondecreasing,
in other words, $\Delta_n$ with $[0,1]$ replaced by $[a,b]$;
denote $\Delta_n^{(L)} := \Delta_n^{[-L,L]}$.
Let $\Delta^{[a,b]}$ denote the disjoint union of the interiors,
also topologized by attachment.  Each $\Delta^{[a,b]}$ embeds 
naturally in $\Delta^{[c,d]}$ for $[a,b] \subseteq [c,d]$.  There 
are natural projections $\pi_{[c,d],[a,b]} : \Delta^{[c,d]} \to 
\Delta^{[a,b]}$ that
ignore points outside of $[a,b]$ and act as one-sided inverses
to the natural embeddings.  We denote by $\pi_{[a,b]}$ the projection
from the inverse limit to $\Delta^{[a,b]}$.  Because $[a,b]$ may now
vary, we will sometimes need the notation $\Phi_{t,[a,b]} (\xx)$ 
instead of $\Phi_t (\xx)$ for the time $t$ map on $\Delta^{[a,b]}$.

Let $\CM$ denote the inverse limit of these projections as 
$M \to \infty$.  Thus, for us, an infinite point configuration 
$\eta$ is a compatible collection of its finite projections
rather than a counting measure.  The projection $\pi_{[a,b]}$ induces 
a projection from $\CM$ to trajectories on $\Delta^{[a,b]}$ via 
$(\pi_{[a,b]} \omega) (t) := \pi_{[a,b]} (\omega (t))$.  Particles 
in the configuration $\pi_{[a,b]} \omega (t)$ die as they pass out of 
$[a,b]$ and are born as they pass into $[a,b]$; this does 
not create a discontinuity in $\pi_{[a,b]} \omega$ because of the 
identifications $\iota_{n,0}$ and $\iota_{n,n+1}$ which allow 
creation and destruction at the endpoints of $[a,b]$.  
In this context, the Poisson measure with intensity $\lambda$,
denoted $\pois(\lambda)$, is the probability measure on $\CM$
that projects under each $\pi_{[a,b]}$ to $\pois (\lambda , [a,b])$.

Let $\TM$ denote the space of continuous trajectories on $\CM$,
that is, continuous maps $\omega : \R^+ \to \CM$.  Topologize
$\TM$ by uniform convergence on compact time intervals $[0,T]$.
Here, by the inverse limit construction, uniform convergence means
uniform convergence of each projection $\pi_L \omega (\cdot)$
on each compact time interval $[0,T]$.
We care only about a tiny subset of $\TM$, namely those trajectories
that might arise as limits of trajectories $\Phi_{t,L} (\xx)$.
To avoid having to craft arguments for paths that might have
exotic behaviors, we define this closure, $\TMm$.  Formally,
for each $n$ and $L$ and each point $\xx \in \Delta_n^{(L)}$,
the trajectory $\{ \Phi_{t,L} (\xx) \}$ lifts to an element of 
$\TM$, the one in $\pi_L^{-1}$ with no points outside of $[-L,L]$.
Let $\TMo$ denote the set of such trajectories and 
let $\TMm$ denote the closure of $\TMo$ in $\TM$.

Our first aim is to show that elements of $\TMm$ look like trajectories
of $\Phi$: they can never uncoalesce and they satisfy the evolution
rule~\eqref{eq:rule1}.  There is no probability involed here.  A further
goal will be to see that with probability~1 there are no multiple
collisions and that in the limit the number of particles remains
locally bounded in probability.  We begin with formal definitions
of collisions, of what it means to solve~\eqref{eq:rule1}, and of
how we count the number of particles in a spacetime rectangle.

\begin{defn}[collision times]
A collision time for a trajectory $\omega \in \TM$ is a time $t$
such that for some $L,\ee, n$ and $1 \leq j \leq n-1$, the trajectory 
$\pi_L \omega$ is in the interior of $\Delta_n^{(L)}$ for the time 
interval $(t - \ee , t)$ and is in $\iota_{n-1,j} \Delta_{n-1}^{(L)}$ 
for the time interval $[t,t+\ee)$.  The location of the collision at
time $t$ is the $j^{th}$ coordinate of $\omega (t)$.  
Although this will turn out to have probability zero, the definition 
of collision allows for more than one collision at time $t$ or for 
a multiple particle collision.
\end{defn}

\begin{defn}[solution] \label{def:solution}
~~\\[-5ex]
\begin{enumerate}[(i)]
\item A trajectory $\omega \in \TM$ is said to obey~\eqref{eq:rule1} 
if each $\pi_L \omega \in \Delta_n^o$ obeys~\eqref{eq:rule1} at
time $t$ as long as $2 \leq i \leq n-1$ and $t$ is not a coalescing
time for $x_i$.
\item A weak solution with Poisson-$\lambda$ initial conditions is a 
probability measure $\P_\lambda$ on $\Omega$ giving probability 
one to the set of trajectories obeying~\eqref{eq:rule1} and
such that the law $\omega (0)$ under $\P_\lambda$ is 
Poisson with mean $\lambda$.
\item A strong solution with Poisson initial conditions is a map
$\Xi$ from $\CM$ to $\TM$ such that $\Xi (\eta) (0) = \eta$ and
$\Xi (\eta)$ obeys~\eqref{eq:rule1} for $\pois(\lambda)$-almost 
every $\eta$.
\end{enumerate}
\end{defn}

\begin{defn}[occupation]
Let $n(a,t) (\omega) = 1$ if $\omega (t)$ has some coordinate equal
to $a$ and zero otherwise.  Here the trajectory $\omega$ may have 
values in $\CM$ or in any $\Delta^{[a,b]}$, noting that
$n(c,t) (\pi_{[a,b]} \circ \omega) = n(c,t) (\omega)$ when 
$c \in [a,b]$.
Let $n(A,t) := \sum_{x \in A} n(x,t)$ denote the occupation of 
the set $A$ by $\omega (t)$.  Note that for $A = (-L,L)$ the value of 
$n(A,t)$ is the dimension $n$ of the face $\Delta_n^{(L)}$ such that 
$\pi_L (\omega (t)) \in \Delta_n^{(L)}$.
\end{defn}

Most of the remainder of this section is devoted to proving
Theorem~\ref{th:existence}.  This will be proved by taking 
a weak limit of solutions to the finite system.  

\begin{pr} \label{pr:closure}
Paths in $\TMm$ have the following properties.
\begin{enumerate}[(i)]
\item Paths never un-coalesce.  Specifically, fix real $a < b$
and $0 \leq s < t$ and fix $\omega \in \TMm$.  If 
$n(\{ a,b \} , u) = 0$ for all $s \leq u \leq t$, that is, 
no particle enters or exits $[a,b]$ during the time interval $[s,t]$,
then $n([a,b],t) (\omega) \leq n([a,b],s) (\omega)$.
\item Paths in $\TMm$ satisfy~\eqref{eq:rule1}.
\label{item:3}
\end{enumerate}
\end{pr}

\noindent{\sc Proof:}  For~$(i)$,
%
fix $\omega \in \TMm$.  Choose a sequence $\{ \omega_L : L \in \Z^+ \}$
converging to $\omega$ such that $\omega_L$ is a trajectory
$\{ \Phi_{t,L} (\xx_L) \}$ for some point $\xx_L \in \Delta^{(L)}$.  
Fix $s < t$ and an interval $[a,b]$ with $n(\{ a,b \} , u) (\omega) = 0$
for all $u \in [s,t]$.  This implies $n(\{ a,b \} , u) (\omega_L) = 0$
for $L$ sufficiently large.  Denote the initial number of particles
in $[a,b]$ by $r := n([a,b],s) (\omega)$.  We need to show that 
$n([a,b] , t) \leq r$.

\noindent{\bf Claim:} If $\xx \in \Delta^{(L)}$ has coordinates 
$x_i = y < z = x_j$ then the corresponding coordinates of 
$\Phi_{t,L} (\xx)$ differ by at most $(z-y) e^{2t}$.
{\bf Proof:} Up to the first collision time, each gap
$y_k := x_k - x_{k-1}$ increases at rate at most $2 y_k$
due to the velocity of $x_k$ being at most $y_k$
and the velocity of $x_{k-1}$ being at least $-y_k$.
This proves the claim up to the first collision time.
Induction then proves it for all times.

Define the $\ee$-mesh occupation of $[a,b]$ of a point 
$\xx \in \Delta^{[a,b]}$, denoted $n_\ee (\xx)$, to be the 
least $j$ such that the set of coordinates of $\xx$ can 
be covered by $j$ intervals of length no more than $\ee$.  
The fact that $\omega_L \to \omega$ implies that for every 
$\ee > 0$ there exists $M_0 (\ee)$ such that
$$m \geq M_0 \; \Rightarrow \; n_\ee (\omega_L (s)) = r \, .$$
Because the trajectories $\omega_L$ do not cross the endpoints
$\{ a , b \}$ in times in $[s,t]$, each particle at time $t$
comes from one or more coalesced particles at time $s$, hence
from the claim we deduce that
$$n_{e^{2(t-s)} \ee} (\Phi_{t,L} (\xx_L)) \leq n_\ee (\xx_L)  = r \, .$$
Sending $L \to \infty$ and noting that $\pi_{[a,b]} \omega$ is 
the uniform limit of $\pi_{a,b]} \omega_L$ on $[s,t]$, we see that
$n_\delta (\omega (t)) \leq r$ for any $\delta > e^{2(t-s)} \ee$.  
Because $\ee$ is arbitrary, this implies that $n([a,b] , t) \leq r$.

For $(ii)$, let $x_i$ denote a particle
for which $t$ is not a coalescing time. 
Uniform convergence of the $i-1, i$ and $i+1$ components of 
$\Phi_{t,L} (\xx_L)$ implies that both the position 
and derivative of the particle $x_i$ converge uniformly.  
The time derivative in the limiting trajectory is therefore 
the limit of the derivatives for the finite trajectories, 
hence the limit trajectory obeys~\eqref{eq:rule1} at time $t$.  
$\Cox$

\begin{lem} \label{lem:limit}
~~\\[-5ex]
\begin{enumerate}[(i)]
\item Let $\xx^{(L)}$ have distribution $\pois (\lambda , L)$.  
Any weak limit as $L \to \infty$ of the laws of the 
trajectories $\{ \Phi_t (\xx^{(L)} \}$ is a weak solution
with Poisson-$\lambda$ initial conditions.  Consequently if 
the family of pushforwards of $\pois(\lambda , L)$ under
$\xx \mapsto \{ \Phi_{t,L} (\xx) \}$ is tight, then such a
weak solution exists.
\item If for every $t$ the limit $\lim_{L \to \infty} 
\Phi_t (\pi_L (\eta))$ exists $\pois(\lambda)$-almost surely,
then this limit defines a strong solution with Poisson-$\lambda$ 
initial conditions.
\end{enumerate}
\end{lem}

\noindent{\sc Proof:} The weak limit of measures supported on a set $S$
is supported on the closure of $S$.  Therefore, any weak limit of
the trajectories $\Phi_{t,L}$ is in $\TMm$ and hence, by part~$(iii)$
of Proposition~\ref{pr:closure}, obeys~\eqref{eq:rule1}.  The
initial conditions are a weak limit of measures $\pois (\lambda , L)$
which is the law $\pois (\lambda)$ on $\TM$.  This is all that is
needed for~$(i)$ along with the observation that tightness of
this family of laws implies the existence of a limit point.  
Statement~$(ii)$ is an immediate consequence of part~$(iii)$ of
Proposition~\ref{pr:closure} and the definition of a strong solution.
$\Cox$

All the time $t$ marginals of the finite system are explicitly
known.  The only work in establishing tightness is to check that
the supremum over a time interval $s \in [0,t]$ of the number of
particles in a fixed interval $[a,b]$ at time $s$ does not have 
positive mass going to infinity as $L \to \infty$.  We have very 
little information about joint distributions of the process at two 
or more times.  However, we can accomplish what we need by an 
identity that bounds the total occupation of the interval in
terms of particle velocities.  

\begin{pr} \label{pr:total}
Let $N(a,t) (\omega) := \sum_{s \in [0,t]} n(a,s) (\omega)$
denote the cumulative occupation of $\{ a \}$ up to time $t$.
Then for any $n$, any time $t$ and any initial configuration 
$\xx \in \Delta_n^{(L)}$, 
$$\int_{-L}^L N(x,t) \, dx \leq 2 L t \, .$$
\end{pr}

\noindent{\sc Proof:} Fixing $\xx$, the successive collision times
$\tau_1 (\xx), \tau_2 (\xx), \ldots$ may be treated as constants.
Analyzing the trajectories separately on $[0,\tau_1), [\tau_1 , \tau_2)$,
etc., and summing the results, we may assume without loss of 
generality that $t < \tau_1$.  The function $I(a,t)$ is the sum
of $n$ indicators of graphs of functions, these functions being
the trajectories of $x_1 , \ldots , x_n$ as a function of time.
Let $v_j (t) = 2 x_j (t) - x_{j+1} (t) - x_{j-1} (t)$ denote the
velocity of the $j^{th}$ particle at time $t$, where we have used 
$x_j (t)$ to denote the $j^{th}$ coordinate $(\Phi_t (\xx))_j$ of
the system at time $t$.  One has an easy bound
\begin{equation} \label{eq:sum abs}
\sum_{j=1}^n |v_j (t)| \leq 
   \sum_{j=1}^n (x_{j+1} - x_j) + \sum_{j=1}^n |x_j - x_{j+1}|
   \leq 2L \, .
\end{equation}
Let $N(a,t,j) = \# \{ s \in [0,t] : x_j (t) = 1 \}$ count incidences
for just the $j^{th}$ particle.  Banach's Indicatrix 
Theorem~\cite[Th{\'e}or{\`e}mes~1 and~2]{banach} (see also~\cite{natanson})
states that $\int_{-L}^L N(x,t,j) \, dx$ is equal to the total variation
of the trajectory $\{ x_j (s) : 0 \leq s \leq t \}$.  This is just
$\int_0^t |v_j (s)| \, ds$.  Summing gives
$$\int_{-L}^L N(x,t) = \sum_{j=1}^n \int_{-L}^L N(x,t,j) \, dx
   = \sum_{j=1}^n \int_0^t |v_j (s)| \, ds  \, .$$
Everything is nonnegative and we may interchange the integral and sum 
to obtain
$$\int_{-L}^L N(x,t) \, dx = 
   \int_0^t \left ( \sum_{j=1}^n |v_j (s)| \right ) \, ds  
   \leq \int_0^t 2L \, ds = 2Lt \, .$$
$\Cox$

\begin{lem} \label{lem:flux}
For any $n$, any time $t$, any $L$, and any interval 
$[a,b] \subseteq [-L,L]$, 
$$\E \int_a^b N(x,t) \, dx \leq 2 (b-a) t + 2 t \lambda^{-1} e^{2t} \, ,$$
where the expectation is with respect to trajectories whose intial point
is uniform on $\Delta_n^{(L)}$.
\end{lem}

\noindent{\sc Proof:} As in the previous proposition,
\begin{equation} \label{eq:V(s)}
\int_a^b N(x,t) \, dx = \int_0^t V(s) \, ds
\end{equation}
where $V(s) = \sum_j |v_j (s)|$ is the sum of the speeds
of all particles in $[a,b]$ at time $s$.  Instead of 
$V(s) \leq 2L$ for all $s$, we have 
$$V(s) \leq 2(b-a) + (M(s)-b) + (a-m(s))$$
where $M(s)$ is the location of the first particle to the 
right of $b$ at time $s$ and $m(s)$ is the location of the 
first particle to the left of $a$ at time $s$.  Each of
the terms $M(s) - b$ and $a - m(s)$ is an exponential
of mean $\lambda^{-1} e^{2s}$, truncated to $L-b$ and
$a + L$ respectively.  The mean of the truncated exponential
is at most the mean of the untruncated exponential, leading to
$$V(s) \leq 2(b-a) + 2 \lambda^{-1} e^{2s} \, .$$
Plugging this into~\eqref{eq:V(s)} proves the lemma.
$\Cox$

\begin{lem}[velocities are bounded in probability] \label{lem:v bounded}
Fix $[a,b]$ and $T$.  For any trajectory $\omega := \{ \Phi_{t,L} (\xx) \}$ 
in $\TMo$, let $V(\omega)$ denote the maximum absolute velocity
of any particle in $[a,b]$ at any time $t \leq T$.  Then under
Poisson initial conditions, $V$ is bounded in probability, meaning that
$$\pois(\lambda,L) \left ( \{ \omega : V(\omega) \geq M \} 
   \right ) \leq g(M)$$
for some function $g$ going to zero as $M \to \infty$ and depending
on $a,b$ and $T$ but not $L$.
\end{lem}

\noindent{\sc Proof:}
Label each initial particle $x_j$ by the positive real number
$\ell (x_j) := e^{2T} \max \{ x_j - x_{j-1} , x_{j+1} - x_j \}$.
From the claim proved in Proposition~\ref{pr:closure}, we know
that at all times between~0 and~$T$, the gaps between $x_j$ and
$x_{j-1}$ is at most $\ell (x_j)$ provided that we interpret the
gap to be zero if the two particles have coalesced.  The same
is true of the gap between $x_j$ and $x_{j+1}$.  It follows
that the two gaps adjacent to a particle $x$ at any time up to $T$
is at most the maximum, call it $W(x)$, of $\ell (x_j)$ over all 
particles $x_j$ that have coalesced into $x$.  Because the velocity 
of any particle at any time is bounded by the maximum adjacent gap, 
the position of any particle $x$ at time $t \leq T$ must be within 
distance $t W(x)$ of the initial position of some $x_j$ that has
coalesced into $x$.  If the particle $x$ at time $t \leq T$ is in
the interval $[a,b]$, and if $x_j$ is the initial position of 
a particle that has coalesced into $x$, then the distance 
$d(x_j , [a,b])$ from $x_j$ to the interval $[a,b]$ must be 
at most $T W(x)$.  Choosing $j$ for which $W(x) = \ell (x_j)$,
we see that $d(x_j , [a,b]) \leq T \ell (x_j)$.  Therefore,
\begin{equation} \label{eq:sup ell}
V(\omega) \; \leq \; V_* (\omega) \; := \;  
   T \cdot \sup \{ \ell (x_j) : d(x_j , [a,b]) \leq T \ell (x_j) \} \, .
\end{equation}

Couple the distributions of $V_*$ under the laws $\pois (\lambda , L)$
as $L$ varies, by taking the restrictions to $[-L,L]$ of a single
pick $\eta$ from the Poisson measure $\pois(\lambda)$ on $\R$.  
The function $\ell (x_j)$ for particles of $\pi_L \eta$ are at
most what they are for the corresponding particles of $\eta$,
with equality except for the first and last particles in $[-L,L]$.
The conclusion of the lemma will therefore follows once we show that
for infinite configurations, 
$$\pois (\lambda) \left [ \{ \eta : V_* (\eta) = \infty \} \right ] = 0 \, .$$
This follows from the fact that the supremum in~\eqref{eq:sup ell} 
is taken over an almost surely finite set.  This, in turn, is a 
simple consequence of the Borel-Cantelli lemma, once one observes 
that the probabilities $q_n$ and $q_n'$ are summable, where
$q_n$ is the probability of existence of an initial particle 
$x$ in $[b+n , b+n+1]$ with $T \ell (x) \geq n$ and $q_n'$ is
the analogous probability for a particle in $[a-n-1,a-n]$.
$\Cox$

\noindent{\sc Proof of Theorem}~\ref{th:existence}:
Fix $L , \lambda > 0$.  On the probability space of trajectories
with initial distribution $\pois (\lambda , L)$, define random
variables $Y = Y(L,a,b,t)$ to be the maximum number of particles
in the interval $[a,b]$ at any time $s \in [0,t]$.  Any particle
contributing to this number is either in $[a-1,b+1]$ at time zero
or crosses the interval $[a-1,a]$ or crosses the interval $[b,b+1]$.
A particle crossing an interval $J$ of length~1 contributes at least~1
to $\int_J N(x,t) \, dx$.  Therefore, 
$$\P (Y > 3y) \leq \P (Y_0 > y) + \P(Y_1 > y) + \P (Y_2 > y)$$
where $Y_0$ is the initial number of particles in $[a,b]$,
$Y_1 = \int_{a-1}^a N(x,t) \, dt$ and $Y_2 = \int_b^{b+1} N(x,t) \, dt$. 
Each of $Y_0, Y_1$ and $Y_2$ is bounded in expectation by a constant
depending on $a,b,t$ and $\lambda$ but not $L$.  Therefore,
$\P (Y > 3y) \leq C / y$, finishing the proof of tightness.

Tightness of these variables for all $a,b$ is tightness
of the pushforwards in part~$(i)$ of Lemma~\ref{lem:limit},
the lemma then implying existence of the weak solution which
we denote henceforth $\{ \Psi_t \}$.  
By construction the two-sided derivative of each particle's position 
exists away from its collision times and~\eqref{eq:rule1} is satisfied.

To argue for almost sure distinctness of collision times, fix 
$\lambda$ and $T$ and the interval $[a,b]$.  For $\ee > 0$, 
define the event $G_\ee$
of an $\ee$-almost multiple collision by time $T$ on $[a,b]$
to be the event that for some times $s$ and $t \in [s,s+\ee]$
there are distinct pairs of particles whose positions differ by 
at most $\ee$ (the two pairs can share one particle but not both).
If $\omega_L \in \TMo$ with $\omega_L \to \omega$, and $\omega$
has a multiple collision at some time $t \leq T$, then for all $\ee > 0$
there will be $L_0 (\ee)$ such that $L \geq L_0$ implies
$\omega_L$ has an $\ee$-almost multiple collision.  Therefore, 
to finish the proof of the theorem, it suffices to show that
the probability of $G_\ee$ under any $\pois (\lambda , L)$ 
is bounded above by some function $\kappa (\ee)$ going to 
zero as $\ee \to 0$, independent of $L$.

Fix an integer $M > 0$ and let $H_M$ denote the event that the maximum
absolute velocity of any particle in $[a,b]$ up to time $T$
is at most $M$.  Let $S$ be the set of multiples of $\ee$
in $[0,T]$.  Then, on $H_M$, the event of an $\ee$-almost multiple 
collision implies that for some time $s \in S$, there are two distinct
pairs of particles (possibly sharing one particle) within 
distance $(M+1) \ee$ at time $s$.  Any pair of points in 
$[a,b]$ within distance $(M+1)\ee$ of each other is in the
same interval $[x , x+2(M+1)\ee]$ for some $x = k (M+1) \ee$,
$0 \leq k \leq 2(b-a) / ((M+1)\ee)$.
Therefore, on $H_M$, the event of an $\ee$-almost collision is contained
in the union of at most $|S| (2 (b-a) \ee^{-1} / (M+1))^2$ events that
two distinct pairs of points at a specific time $s$ are
in two distinct intervals $[x_1 , x_1 + 2 (M+1) \ee]$ and
$[x_2 , x_2 + 2 (M+1) \ee]$, together with the $|S| 2(b-a) \ee^{-1} / (M+1)$
events that two distinct pairs of points at a specific time $s$
are together in the same interval $[x,x+2(M+1)\ee]$.

We bound this from above using the fact that the time $s$ marginal 
of the number of points in any interval of length $c$ is a 
Poisson of mean $\lambda e^{-2s} c \leq \lambda c$.  The probability
of a Poisson with mean $\nu$ being at least~2 is at most $\nu^2 / 2$
and the probability of it being at least~3 is at most $\nu^3 / 6$.
Therefore, applying the bound with $c = 2 (M+1) \ee$, and using
the bound $|S| \leq T \ee^{-1} + 1$, we see that the probability 
of an $\ee$-almost collision is at most 
$$(T \ee^{-1} + 1) \frac{2(b-a)}{M+1} \ee^{-1} \frac{(2 (M+1) \ee)^3}{6}
   + (T \ee^{-1} + 1) \left ( \frac{2(b-a)}{M+1} \ee^{-1} \right )^2
   \left ( \frac{(2 (M+1) \ee)^2}{2} \right )^2 \, .$$
Clearing away irrelevant stuff, we see that
$$\pois (\lambda, L) (H_M \cap G_\ee) \leq  C T (b-a)^2 M^2 \ee \, .$$
By Lemma~\ref{lem:v bounded}, the probability of $H_M^c$ is bounded 
above by $g(M)$.  Choosing $M = \ee^{-1/3}$ we see that
$$\pois (\lambda, L) ( G_\ee) \leq  C T (b-a)^2 \ee^{1/3} 
   + g(\ee^{-1/3}) \, .$$
This goes to zero as $\ee \downarrow 0$, finishing the proof
of Theorem~\ref{th:existence}.
$\Cox$

\section{Further comments and questions}

The most glaring absence of a result concerns strong solutions.
The following conjecture, together with Lemma~\ref{lem:limit}
would imply the existence of a strong solution to~\eqref{eq:rule1}.
A proof seems not too far off via results along the lines of 
Lemma~\ref{lem:flux}.

\begin{conj}[strong solution] \label{conj:converge}
For $\P$-almost every $N \in \Omega$, the limit
$$\Psi_t (N) := \lim_{L \to \infty} \Psi_t^{(L)} (N)$$
exists and defines a trajectory $t \mapsto \Psi_t (N)$.
\end{conj}

Conjecture~\ref{conj:converge} also implies the following 
restatement of Theorem~\ref{th:poisson invariant} in terms
of the evolution $\{ \Psi_t \}$.  
\begin{equation} \label{eq:push forward}
\pois (\lambda) \circ \Psi_t^{-1} = \pois (\lambda e^{-2t}) \, .
\end{equation}

Numerical simulation suggests a substantially more general result than 
those obtained in this paper.  We believe that the Poisson distribution
is an attracting fixed point for the dynamics: any reasonable initial 
measure, if allowed to evolve under these dynamics and rescaled, should
converge Poisson.  There is some numerical evidence for this.
In one experiment, 250 million points were initially placed on the 
unit interval with uniform density and periodic boundary conditions.  
New points were then placed 
halfway between pairs of adjacent points, and the initial points 
were removed.  The resulting configuration is known as the 1D 
Poisson-Voronoi configuration.  Its cell sizes are one-dependent
but not independent, with marginal density $4x \exp (-2x)$.  
Figure~\ref{fig:numerical} illustrates histograms of the normalized
cell sizes as time grows.  As increasingly many of the cells disappear, 
the pdf approaches the exponential distribution, the cell size
marginal for the Poisson.  The data is also compatible with 
asymptotic independence of neighboring cell sizes.  Additional 
supporting data from this initial configuration and several other
initial configurations are reported in~\cite{lazar2011evolution}.

\begin{center}
\includegraphics[height=3in]{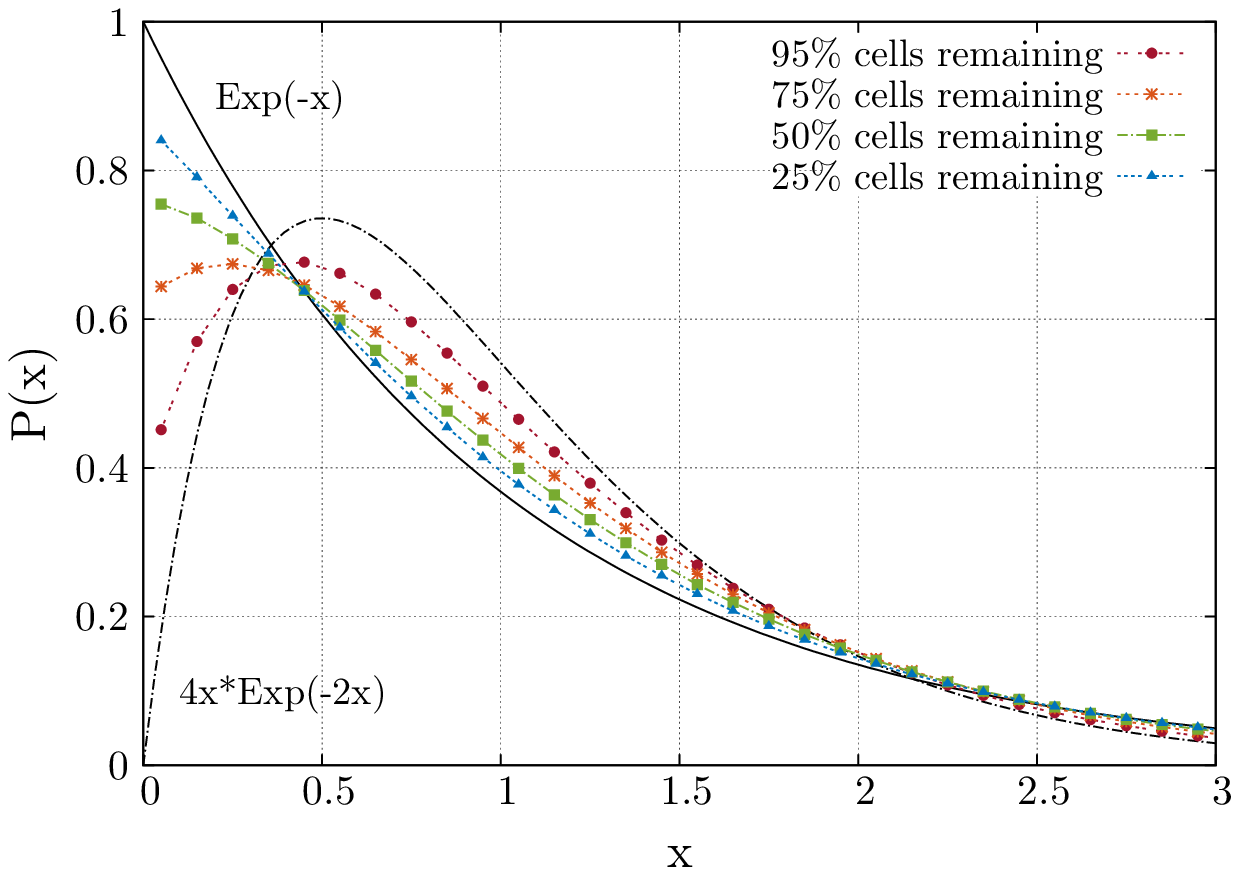}
\label{fig:numerical}
\end{center}

To make this into a precise conjecture,
let $\nu$ be a probability measure on $\R$ with finite mean and
let $Q_\nu$ be the stationary renewal process on $\R$ with renewal
distribution~$\nu$.  Intuitively, $Q_\nu$ is the law of a random
configuration of points whose gaps are IID~$\nu$.

\begin{conj} \label{conj:attractor}
For any $\nu$ on $\R$ with mean $m < \infty$, there is a weak solution 
$\P_\nu$ to~\eqref{eq:rule1} with initial conditions $Q_\nu$.
Let $\nu_t$ be the time $t$ law of such a weak solution,
with space rescaled by $e^{-2t}$.  Then $\nu_t \to \pois (1/m)$
weakly as $t \to \infty$.
\end{conj}

This conjecture is the concrete assertion of the universality 
mentioned in the abstract.  
While a proof seems farther off for this conjecture than for 
Conjecture~\ref{conj:converge}, we hope to attack it via 
the dual process introduced in Proposition~\ref{pr:reversals}.  
While the forward evolution is deterministic, the time-reversed
process has randomness.  We believe this will help us to show
that running the dual process back in time $t$ units produces, as 
$t \to \infty$, a configuration asymptotically independent from
any $Q_\nu$.  

\section*{Acknowledgements}

The first author would like to thank Benjamin
Matschke and Amanda Redlich for helpful discussions.
The second author would like to thank Ryan Hynd and 
Herman Gluck for helpful conversations.  Thanks also to 
Moe Hirsch and Joel Robbin for providing references to the 
works of J.\ Aubin and P.\ Hartman. \\[-9ex]

\bibliographystyle{alpha}
\bibliography{RP}

\end{document}